\documentclass[default]{sn-jnl}
\usepackage{amstext,amsmath,amsfonts,amssymb,amsbsy,bm,color,graphicx,subfigure,natbib}
\usepackage{enumitem}
\jyear{2022}%
\theoremstyle{thmstyleone}%
\newtheorem{theorem}{Theorem}
\newtheorem{proposition}[theorem]{Proposition}%

\theoremstyle{thmstyletwo}%
\newtheorem{example}{Example}%
\newtheorem{remark}{Remark}%

\theoremstyle{thmstylethree}%
\newtheorem{definition}{Definition}%

\def\R{\mathbb{R}}

\begin{document}

\title[Time discretization in the solution of parabolic PDEs with ANNs]{The effect of time discretization on the solution of parabolic PDEs with ANNs
}


\author[1]{\fnm{Francesco} \sur{Calabr\`o}}\email{francesco.calabro@unina.it}

\author[1]{\fnm{Salvatore} \sur{Cuomo}}\email{salvatore.cuomo@unina.it}

\author*[1]{\fnm{Daniela} \sur{di Serafino}}\email{daniela.diserafino@unina.it}

\author[1]{\fnm{Giuseppe} \sur{Izzo}}\email{giuseppe.izzo@unina.it}

\author[1]{\fnm{Eleonora} \sur{Messina}}\email{eleonora.messina@unina.it}

\affil[1]{\orgdiv{Dipartimento di Matematica e Applicazioni ``R.~Caccioppoli''}, \orgname{Università degli Studi di Napoli Federico II}, \orgaddress{\street{Monte Sant'Angelo, via Cintia}, \city{Napoli}, \postcode{80126}, \country{Italy}}, Member of the INdAM-GNCS Research Group}


\abstract{
We investigate the resolution of parabolic PDEs via Extreme Learning Machine (ELMs) Neural Networks, which have a single hidden layer and can be trained at a modest computational cost as compared with Deep Learning Neural Networks. Our approach addresses the time evolution by applying classical ODEs techniques and uses ELM-based collocation for solving the resulting stationary elliptic problems. In this framework, the $\theta$-method and Backward Difference Formulae (BDF) techniques are investigated on some linear parabolic PDEs that are challeging problems for the stability and accuracy properties of the methods.
The results of numerical experiments confirm that ELM-based solution techniques combined with BDF methods can provide high-accuracy 
solutions of parabolic PDEs.
}




\maketitle

\section{Introduction}

Mesh-based schemes are widely discussed numerical approaches for solving Partial Differential Equations (PDEs). Finite Difference Methods (FDMs), Finite Element Methods (FEMs), and Finite Volume Methods (FVMs) all belong to this class of numerical methods \cite{LiborAQ}. They require the generation of grid points or elements, usually the discretization of differential/integral formulas, and the solution of discrete equations, often with iterative algorithms. Mesh-based approaches suffer from several problems: firstly, the shape complexity of the computational domain where the grid generation itself could become very difficult or even infeasible; moreover, the problem discretization could create a bias between the mathematical nature of the PDE and its approximating model.

In recent years, starting from \cite{lagaris1998artificial}, Artificial Neural Networks (ANNs) have been considered an interesting alternative methodology to overcome the drawbacks of mesh-based numerical schemes. ANNs are adopted as basis functions to compute solutions of PDEs, i.e.~the approximate problem solutions are determined by a learning approach that consists in defining an optimization algorithm in which losses due to ANN approximations of PDEs and boundary conditions (BCs) are minimized. This approach requires sampling points inside the domain and on the boundary, which can be randomly selected. When the numerical approximation of the unknown solution is guided by the resolution of the underlying PDE, the methods are referred to as physics informed \cite{karniadakis2021physics,AA}.

Recently, it has been demonstrated that PDEs can be solved by considering a specific ANN called Extreme Learning Machine (ELM), see \cite{huang2015trends}. An ELM is a feed-forward neural network with a single hidden layer that randomly assigns the input layer weights and analytically determines the output weights. Thanks to this architecture, the weights of the hidden layer need not be learned. This makes ELMs faster than typical deep neural networks, where optimization methods may lead to prohibitively slow learning speeds. We point out that ELMs are variants of the random projection networks originally proposed in~\cite{RPNN}, and a discussion on the relation beween this and other theories regarding random networks can be found in the review paper~\cite{Rev_RPNN}. Overall, randomized Neural Networks boost the learning task with benefits on the numerical scheme in terms of efficiency, while maintaining high accuracy. ELMs have been successfully applied for solving Ordinary Differential Equations (ODEs) and stationary PDEs \cite{Cal_1,Dwi2020_PhysicsInformedExtreme_SriDS,ZZ2,Cal_2,schiassi2021extreme}. ELMs are mesh-free methods and thus they can easily address complex geometries of PDE domains~\cite{Dwi2020_PhysicsInformedExtreme_SriDS}. Moreover, ELMs are universal approximators \cite{huang2006universal,huang2006extreme} and hence can potentially represent any continuous function.

In this paper we consider ELMs for the numerical resolution of a linear parabolic equation. Following what is discussed in~\cite{Cal_1}, we use collocation for the resolution of stationary elliptic problems, but address separately the time marching. The resolution of the elliptic problem has provided a very good accuracy with a modest computational cost in situations where classical methods fail to give good results. Moreover, it has very nice properties of generalization and has been applied to nonlinear problems in~\cite{Cal_2}. As usually done in time-dependent problems, we aim to separate the dependence on the time by representing the unknown solution as a combination of stationary functions where the coefficients of the combination change in time.

Now we show how collocation can be applied to our differential problem. Let $u(t,x)$ be a solution to the following:
\begin{equation}\label{eq:1}
    \frac{\partial u}{\partial t} = \mathcal{L}(u) + f(t,x), \quad \mathcal{B}(u) = g(t,x), \quad
\end{equation}
where $\mathcal{L}$ is intended as a linear elliptic operator acting only on the spatial variables, $\mathcal{B}$ is a boundary operator and $f,g$ are given source and boundary data. We are interested in the numerical resolution of equation~\eqref{eq:1}, in particular in a solution that can be written as an ANN of ELM-type in space: \begin{equation}\label{eq:elm_sol}
u(t,x)\approx u^h(t,x)=\sum_i w_i(t)\sigma_i(x).
\end{equation}
The main assumption of our collocation method is that we are able to compute
\begin{equation}\label{eq:hyp_coll}
\mathcal{L}(\sigma_i(x)),\;\;\mathcal{B}(\sigma_i(x)).
\end{equation}
We denote by $\{x_j\}$ the set of points where the previous quantities are evaluated. Then, two strategies are possible for the time marching, and these lead to different classes of methods, as specified next.
\begin{enumerate}
\item
In the first class we focus on the unknown functions $w_i(t)$ after applying the operator $\mathcal{L}$ (or some approximation to it) to the functions $\sigma_i$. In this case we make a semi-discretization in space and then solve a system of ODEs, which is linear if the PDE is linear.
In order to describe this approach, we apply the operators $\frac{\partial}{\partial t}$ and $\mathcal{L}$ to the solution in~\eqref{eq:elm_sol}:
\begin{equation}\label{eq:4}
\frac{\partial}{\partial t}u^h(t,x)=\frac{\partial}{\partial t}\left(\sum_i w_i(t)\sigma_i(x)\right)= \sum_i \sigma_i(x) \frac{d}{d t}w_i(t),
\end{equation}
\begin{equation}\label{eq:5}
\mathcal{L}(u^h(t,x))=\mathcal{L}\left(\sum_i w_i(t)\sigma_i(x)\right)= \sum_i w_i(t)\mathcal{L}(\sigma_i(x)).
\end{equation}
Now we impose that \eqref{eq:4} and \eqref{eq:5} solve exactly problem~\eqref{eq:1} at the given points~$x_j$. Let $G = (g_{ij})$ be the Gram matrix, where $g_{ij}=\sigma_i(x_j)$, and let $L=(l_{ij})$ be the matrix where $l_{ij}=\mathcal{L}(\sigma_i(x_j))$. With these positions, the final resolution corresponds to solving
\begin{equation}
G \frac{d}{dt} w(t) = L w(t) + f(t)
\end{equation}
where $w(t)=(w_i(t))_i$ is the vector of the unknowns and $f(t) = (f(t,x_j))_j$ is the vector of the source terms. We refer to~\cite{thomeeCapitolo} for details.
\item
In the second class of methods we define
\begin{equation} \label{eq:sol_meth_2}
\tilde{u}^{[n]}(x) = \sum_i w^{[n]}_i \sigma_i(x),
\end{equation}
where $w^{[n]}_i = w_i(t_n)$ and $\{t_n\} $ is a discretization of the time interval with fixed time spacing $\Delta t$. Then the unknowns are the coefficients $w^{[n]}_i$, which are computed by solving a boundary problem. In this case we first perform a semi-discretization in time and then solve an elliptic stationary problem, which is linear if the PDE is linear.
As an example, let us consider the well-known $\theta$-method, which includes the explicit Euler method for $\theta=0$, the backward (implicit) Euler method for $\theta=1$ and the trapeziodal (Crank-Nickolson) method for $\theta=1/2$. At each time step, we look for the function $\tilde{u}^{[n]} (x)$ in the hypothesis that $\tilde{u}^{[n-1]} (x)$ is known, i.e.~given or computed. Then, the method solves the following:
\begin{equation}\label{theta}
    \frac{\tilde{u}^{[n]} (x)-\tilde{u}^{[n-1]} (x)}{\Delta t} = \theta \left[\mathcal{L}(\tilde{u}^{[n]}) + f(t_n,x)\right]  + (1-\theta)\left[ \mathcal{L}(\tilde{u}^{[n-1]}) + f(t_{n-1},x) \right],
\end{equation}
that is
\begin{align*}
    - \theta \Delta t \mathcal{L}(\tilde{u}^{[n]}) + & \tilde{u}^{[n]} (x) =  \\
    & \tilde{u}^{[n-1]} (x) + \theta \Delta t  f(t_n,x)  + (1-\theta) \Delta t \left[ \mathcal{L}(\tilde{u}^{[n-1]}) + f(t_{n-1},x) \right].
\end{align*}
The last equation can be regarded as an elliptic equation.\footnote{Notice that it reveals to an explicit method -- the forward Euler method -- in the case $\theta=0$.} The right-hand side consists of given functions and quantities that can be explicitly evaluated.
\end{enumerate}

The first class is more natural for problems where the diffusion is dominated by the transport, leading to hyperbolic-type behaviors, or where stiffness arises in time. In these cases one can easily adopt different time discretizations. Moreover, the first class of methods can be profitably applied when the matrix involved in the problem has a structure that can be exploited. The introduction and study of different time-marching procedures has been proposed also for the parabolic case, for both traditional methods~\cite{AA1_gen,AA2_gen} and newer ones, such as Isogeometric Analysis~\cite{AA1_iga,AA2_iga,AA3_iga} and Physics-Informed Neural Networks~\cite{AA}.

On the other hand, the differential problem that has to be solved in the case of ELMs involves a usually dense (and not a-priori banded) matrix and the properties of this matrix are difficult to establish. In the authors' opinion, the ELM collocation can be more succesfully applied in the second class of time-marching methods, as we discuss in this paper. A main advantage is that one can apply different strategies for the resolution of the elliptic problem at different time steps, e.g.~by changing the collocation points. 

\bmhead{Our contribution}

The aim of this work is to explore how the time approximation affects the resolution of time-dependent parabolic PDEs when the time marching is made according to the second class of methods described before.
We start from the above-mentioned $\theta$-method and then apply Backward Differentiation Formula (BDF) methods to overcome some difficulties, see~\cite{Cash,TR_BDF2}. BDF methods show good accuracy and nice convergence properties while keeping the computational cost unchanged with respect to the implicit $\theta$-method. Moreover, when the time step is chosen to be constant, they are easy to implement.

The numerical tests fully confirm the reported results concerning stability and order of convergence of the proposed time-discretization schemes also  when combined with the ELM-based collocation method. 
In the case of the parabolic problems presented here, this reveals that BDF methods are to be preferred to the well-established Backward Euler or Trapezoidal rules because the latter are unable to deal with accuracy requirements, as they are slowly convergent and can suffer from order reduction or numerical instabilities.

\bmhead{Structure of the paper and notation}

The rest of this paper is organized as follows. In Section~\ref{sec:ann-approximators} we briefly discuss the use of Single hidden-Layer Feedforward Networks (SLFNs) as function approximators, focusing on the selection of the Activation Functions (AFs) and on the application of ELMs as collocation methods for the solution of PDEs. In Section~\ref{sec:time-marching} we present BDF solvers, which are used as time-marching methods in the ELM-based resolution. The results of numerical experiments, reported in Section~\ref{sec:experiments}, show that our approach is able to provide accurate solutions to linear, but challenging, parabolic PDE problems, according with the order of the time-marching procedure and the theoretical results on ELMs. Some conclusions are given in Section~\ref{sec:conclusions}.

In the following, vectors are written in boldface, i.e.~$\bm{v}$, while scalars are written in lightface, i.e.~$v$. The $i$-th entry of a vector $\bm{v}$ is denoted $v_i$ and the scalar product of $\bm{v}$ and $\bm{w}$ is denoted $\bm{v} \cdot \bm{w}$. Finally, $\| \cdot \|$ indicates either the Euclidean norm of a vector or the $L^2$-norm of a function.

\section{Use of ANNs as approximators\label{sec:ann-approximators}}

Today many researchers agree that ANNs learn to approximate functions. More in detail, ANNs are  techniques for estimating an unknown function using available observations or collocation points from the function domain~\cite{higham2019deep,pinkus}. The function to be estimated, which transforms inputs to outputs, is often referred to as the target function.

Here we discuss the application of SLFNs with random hidden nodes to a differential problem. A SLFN with $\bar{N}$ hidden nodes and AFs $\phi: \R \to \R$ is a function $\mathcal{F}:\R^d\to\R^m$ defined as follows:
\begin{equation}\label{ANN_vec}
\mathcal{F} (\bm{x}) = \sum_{i=1}^{\bar{N}} \bm{w}_i \phi_i(\bm{x}) = \sum_{i=1}^{\bar{N}} \bm{w}_i \phi(\bm{\alpha}_i \cdot \bm{x} + \beta_i), 
\end{equation}
\noindent
where $\bm{\alpha}_i = (\alpha_{i1}, \alpha_{i2}, \dots, \alpha_{id})$ is the weight vector linking the input nodes to the $i$-th hidden node, $\bm{w}_i=(w_{i1}, w_{i2}, \dots, w_{im})$ is the weight vector linking the $i$-th hidden node to the output nodes, and $\beta_i$ is a bias of the $i$-th hidden node. Our aim is to study whether the SLFN $\mathcal{F}(\bm{x})$ fits the data.
More formally, given $M$ arbitrary  couples $(\bm{x}_j,\bm{\tau}_j)$, where
$$
\bm{x}_j=(x_{j1}, x_{j2}, \dots, x_{jd})\in \R^d, \quad \bm{\tau}_j=(\tau_{j1}, \tau_{j2}, \dots, \tau_{jm}) \in \R^m,
$$
and $\bm{x}_j \ne \bm{x}_i$ for $j \ne i$, these are $M$ samples where the SLFN passes with zero error:
$$
\mathcal{F}(\bm{x}_j) = \bm{\tau}_j, \quad j = 1, 2, \dots ,M .
$$

\noindent
SLFNs, as a special case of ANNs, are universal approximators in the sense that a feedforward network with a linear output layer and at least one hidden layer with any nonlinear activation function can approximate with good accuracy a given function from a finite-dimensional space to another, provided that the network has enough hidden nodes, see, e.g., \cite{cybenko1989,hornik1989} and the recent review~\cite{kratsios2021universal}. Here our main interest is on scalar functions, corresponding to $m=1$. In this case, the vectors $\bm{\tau}_i$ become the scalars $\tau_i$ and function~\eqref{ANN_vec} reads
\begin{equation} \label{ANN}
\mathcal{F} (\bm{x}) = \sum_{i=1}^{\bar{N}} w_i \phi_i(\bm{x}) = \sum_{i=1}^{\bar{N}} w_i \phi(\bm{\alpha}_i \cdot \bm{x} + \beta_i), 
\end{equation}
where the weights $w_i$ are scalar too.

For the sake of completeness, we report next a version of the universal approximation theorem for the SLFN~\eqref{ANN}, see~\cite[Theorem 3.1]{pinkus}.
\begin{theorem}\label{th:univ_approx}
Let $\mathcal{F}$ be a SLFN function as in \eqref{ANN}, where $\phi$ is not a polynomial. 
For any continuous function $f: \R^d \to \R$, any compact subset $K\subset \R^d$ and any $\varepsilon>0$, there exist weights and biases $\bm{w}_i, \bm{\alpha}_i, \beta_i$ such that
\begin{equation} \label{eq:approx}
\max_{\bm{x}\in K}  \mid \mathcal{F}(\bm{x})-f(\bm{x}) \mid < \varepsilon .
\end{equation}
\end{theorem}

\noindent
This result states that a single-layer network is enough to have universal approximation. In order to get a good approximation, we attempt to minimize the left-hand side of \eqref{eq:approx} or, more generally, a measure of the distance between $\mathcal{F}$ and $f$, usually referred to as \emph{loss function}.

Additional hidden layers can help model intricate nonlinear dynamics \cite{Sen2020_ReviewDeepLearning_BasSBS,mhaskar2016deep}.
In deep networks with many layers, an important role is played by Physics Informed Neural Networks (PINNs), often adopted for solving real problems. However, deep learning approaches are characterized by high training costs and by efficiency issues, while we use shallow networks with random projection neurons, leading to ELM networks.

\subsection{Selection of the activation functions\label{sec:selection}}

A crucial aspect of the learning approach is the selection of AFs, since they can significantly affect the accuracy and efficiency of an ANN. Commonly used AFs are ReLU, leaky ReLu, Sigmoid, Tanh (see, e.g., \cite{GF,higham2019deep}), but other choices have been also discussed, e.g., in~\cite{Sun2020_SurrogateModelingFluid_GaoSGPW,He2020_PhysicsInformedNeural_BarHBTT,Che2021_DeepLearningMethod_ZhaCZ}.

In many problems, it is necessary to rescale the PDE to a dimensionless form. In this case, for the selection of an AF, it is recommended to pick a fixed range, such as $[0,1]^d$, rather than to consider the whole domain where the problem is defined.
Moreover, the regularity of an ANNs can be obtained by utilizing smooth activation functions like the hyperbolic tangent or a sigmoid. However, for all non-polynomial AFs an interpolation result can be stated~\cite[Theorem 5.1]{pinkus}.

\begin{theorem}\label{th:interp}
Let $\mathcal{F}$ be a SLFN function as in \eqref{ANN}, where $\phi$ is not a polynomial. 
For any $M$ distinct points $\bm{x}_j$ and associated data $\tau_j$, there exists a choice of $M$ coefficients $\bm{\alpha}_i=(\alpha_{i1},\dots,\alpha_{id} )\in \R^{d}, {w}_i\in \R, {\beta}_i\in \R$ such that
\begin{equation*}
\sum_{i=1}^{{M}} {{w}_i} \phi( \bm{\alpha}_i \cdot \bm{x}_j+\beta_i)=\tau_j, \quad j=1,\dots,M . 
\end{equation*}
\end{theorem}
\noindent
In other words, for an SLFN architecture with $N=M$ hidden nodes, it is possible to approximate $M$ samples with zero mean error.



In this work we choose the AFs in the class of sigmoid functions. A sigmoid usually takes values between 0 and 1 and it is widely used for models where a probability has to be predicted as an output. In our case, sigmoids are a good choice on a collocation basis as they are differentiable. This means we can find the slope of the sigmoid curve. Moreover, this function gives an interesting advantage in terms of computational time for the training phase of the neural network. 

Specifically, our choice of the AF is the logistic sigmoid function:
\begin{equation}
\sigma_i(\bm{x})=\sigma(\bm{\alpha}_i \cdot \bm{x} + \beta_i)= \frac{1}{1+\text{exp}(-\bm{\alpha}_i \cdot \bm{x}-\beta_i)} .
\label{def:sigmai}
\end{equation}
The derivatives of \eqref{def:sigmai} with respect to the independent variable $\bm{x}$ can be easily computed.
Note that if one takes two functions $\sigma_i$ and $\sigma_j$, where at least one of the parameters is different, then these functions are linearly independent, see \cite{ito1996nonlinearity}. Moreover, each $\sigma_i$ is a planar-wave Ridge function~\cite{pinkus}, so that the behavior of the function can be derived easily as an extension of the univariate case $d=1$. With this simplification, we can state the following:
\begin{itemize}
\item $\sigma_i$ has an inflection point at $x=-\dfrac{\beta_i}{\alpha_i}$, which we call the \textit{center} $C_i$ of the sigmoid function;
\item $\sigma_i$ is monotone, $\displaystyle \lim_{x\to -\infty}=0$ and $\lim_{x\to +\infty}=1$ if $\alpha_i$ is positive, the other way if $\alpha_i$ is negative. Moreover, the range where the values are between 0.05 and 0.95 is $\left[C_i -\frac{2.945}{\alpha_i}, C_i+ \frac{2.945}{\alpha_i}\right]$.
\end{itemize}
Notice that one obtains Heaviside-like functions if the internal weights $\alpha_i$ are large, or
almost-linear functions if the $\alpha_i$'s are small. In our case, the use of both kind of functions can help approximate steep gradients and global behaviors. These functions are an example of AFs in the class of those that verify the hypotheses of Theorems~\ref{th:univ_approx} and~\ref{th:interp}.

\subsection{Shallow networks and ELMs\label{sec:shallow-elms}}

The aim of the overall network is to have nice properties of reproduction while maintaining small the number of unknowns. 
Sparse neural networks have been recently proposed instead of fully-connected architectures to overcome some issues related to learning processes. A Physics-based interpretable sparse neural network architecture for solving PDEs has been analyzed in \cite{Ram2021_SpinnSparsePhysics_RamRR}. It represents a successful tentative to link traditional Deep Neural Networks (DNNs) and meshless methods. The proposed methodology is very efficient in comparison with classical DNNs and it represents a generalized physics-based approach, in the sense that the loss function depends directly on the PDE formulation. Moreover, this model implicitly encodes mesh adaptivity as a part of its training process, leading to novel hybrid algorithms for PDEs. The interpretability of the model is due to a new class of sparse network architectures that generalize traditional meshless methods exactly representing a DNN.

Many authors have suggested investigating shallow networks, e.g.~single hidden-layer networks with an increasing number of neurons in the hidden layer. In the PINNs context, shallow ANNs have been considered a good choice with respect to deep learning methodologies. Among these shallow networks, ELMs~\cite{huang2006extreme} are the ones where internal parameters are fixed randomly and only the external weights $w_i$ are trainable parameters. A fascinating combination of PINNs and ELMs has been investigated in \cite{Dwi2020_PhysicsInformedExtreme_SriDS}, where the authors implemented a model called Physics Informed Extreme Learning Machine (PIELM) for the resolution of stationary and time-dependent linear PDEs. The nice behavior of such networks is related to the universal approximation result, that is valid as for the general SLFNs seen previously.

In \cite[Theorem 2.1]{huang2006extreme} an approximation result needed in our setting is given. We report it next.

\begin{theorem}\label{th:ELM1}
Let $(\bm{x}_i,\tau_i)$, $i=1,\dots, M$, be a set of points such that $\bm{x}_i \ne \bm{x}_j$ if $i \ne j$, and let ${u}_N (\bm{x})= \sum_{i=1}^N w_i \sigma_i(\bm{x}) = \sum_{i=1}^N w_i \sigma(\bm{\alpha}_i \cdot \bm{x} + \beta_i)$ be an ELM network with $N<M$ neurons such that the internal weights $\bm{\alpha}_i$ and the biases $\beta_i$ are randomly generated independently from the data, according to any continuous probability distribution. Then, for all $\varepsilon >0$ there exists a choice of the weights $w_i$ such that
$$
    \| ({u}_N (\bm{x}_i) - \tau_i)_i\| < \varepsilon \quad \mbox{with probability~1},
$$
where $({u}_N (\bm{x}_i) - \tau_i)_i$ denotes the vector with components ${u}_N (\bm{x}_i) - \tau_i$.
Moreover, if $N=M$ then $w_i$, $i=1,\ldots,N$, can be found such that
$$
    \| ({u}_N (\bm{x}_i) - \tau_i)_i \| = 0 \quad \mbox{with probability~1}.
$$
\end{theorem}

\noindent
In particular, the above theorem states that if the number of hidden neurons is equal to the number of data points, then the interpolation error is zero with probability 1. Unlucky cases, i.e.~with probability 0, are related to the unisolvence of the points. In the case $d=1$ the unisolvence hypothesis is included in the request that the points are distinct.

The interpolation property of Theorem~\ref{th:interp} can be extended to a convergence result, as proved in \cite{huang2006universal} (see also \cite[Theorem 2]{huang2015trends}).
\begin{theorem} \label{cor1} 
Let $\phi: \R^d \to \R$ be a continuous function. Then there exist a sequence of ELM network functions ${u}_N (\bm{x}) = \sum_{i=1}^N w_i \sigma(\bm{\alpha}_i \cdot \bm{x} + \beta_i)$ such that: 
\begin{equation*}
 \lim_{N \to \infty} \| \phi - u_N \| = 0.   
\end{equation*}
\end{theorem}

A fundamental challenge in all the works about neural networks is to find an optimal choice of parameters that satisfies a desired tolerance. As commented before, with ELMs it is possible to focus on the optimization of the external weights $w_i$ only, reducing the computational cost and the training time. Thus, the weights $w_i$ can be seen as the coefficients of the linear combination defining ${u}_N (\bm{x})$, while $\bm{\alpha}_i$ and $\beta_i$ are internal weights and biases that yield a variation of the sigmoid function.

Firstly, we present the way in wich the parameters $\bm{\alpha}_i$ and $\beta_i$ are taken, and then we discuss the computation of the weights $w_i$. The internal weights $\bm{\alpha}_i$ are chosen randomly and uniformly in a range that depends on the number of neurons $N$. By fixing the domain of the differential problem to have unitary length and following the analysis in \cite{Cal_1}, we choose
\begin{equation}\label{eq:alpha}
\bm{\alpha}_i = {\text{\textbf{rand}}\left(\left[-\frac{N-10}{10}-4,\frac{N-10}{10}+4\right]\right)} ,
\end{equation}
\noindent
where $\text{\textbf{rand}}([a,b])$ denotes for each $i$ a vector with components sampled from a uniform distribution in the interval $[a,b]$. The biases $\beta_i$ are set so that 
the functions $\sigma_i$ are ``non-flat" in the considered domain, as shown at the end of Section~\ref{sec:selection}. 

The trainable parameters of the network are the weights $w_i$. In general ANNs, the weights are computed by minimizing the loss function, e.g.~by applying stochastic gradient-based approaches that back-propagate the error and adjust the weights through specific directions~\cite{bottou:2018}. More recently, second-order stochastic optimization methods have been widely investigated to get better performances than first-order methods, especially when ill-conditioned problems must be solved, see, e.g., \cite{diserafino2021lsos} and the references therein. Nevertheless, there are still difficulties in using these approaches, such as the setting of the so-called hyperparameters and the significant increase of computing time when the number of nodes in the hidden layer grows.

In our case, by collocating the linear problem~\eqref{eq:1} we obtain a linear dependence on the unknown weights $w_i$. Following Theorems~\ref{th:ELM1} and~\ref{cor1}, we can choose a number of equations smaller than the number of unknowns, obtaining an underdetermined linear system that can be solved as a least squares problem, which plays the role of a loss function. 
By computing the minimum-norm least squares solution~\cite{bjorck1996leastsquares} we not only obtain existence and uniqueness of the weights, but also ensure minimum training error and the smallest norm of the weight vector, which are important properties for an ELM~\cite{huang2006extreme}. This solution can be obtained by using a Complete Orthogonal Decomposition (COD) \cite{bjorck1996leastsquares} of the collocation matrix $C$:
$$
   C P = Q R Z^T,
$$
where $Q$ and $Z$ are orthogonal matrices, $R$ is an upper triangular matrix and $P$ is a permutation matrix. We note that COD has good stability properties even in the case $C$ has (numerical) rank $r < N$~\cite{HoughVavasis1997}. On the other hand, in general COD does not tend to zero the entries of the solution, thus producing functions $u_N$ where most of the coefficients $w_i$ are likely to be nonzero.

Alternative solutions to the undetermined linear system can be computed by using, e.g., the $QR$ factorization with column pivoting~\cite{bjorck1996leastsquares}:
$$
  C P = Q R,
$$
where $Q$ is an orthogonal matrix, $R$ is an upper triangular matrix and $P$ is a permutation matrix. In this case, at least $M-N$ coefficients $w_i$ are set equal to zero, generally yielding a sparser solution than in the previous case.

Preliminary numerical experiments have shown that the difference between the PDE solutions obtained with the two approaches are practically negligible. Therefore, in our tests we compute the minimum-norm least squares solution to the underdetermined linear system.

\section{Time-marching scheme\label{sec:time-marching}}


In the literature on numerical methods for diffusion equations much attention has been paid to the construction and analysis of stable and accurate approximation schemes~\cite{Cash,AA2_gen,Ramos2007}. In this work we compare the trapezoidal~\eqref{theta} and BDF methods according to the second approach described in the introduction, thus we follow the notation introduced in equation~\eqref{eq:sol_meth_2}. Given
\begin{equation}\label{eq:starting}
 \tilde{u}^{[0]}(x),\ \tilde{u}^{[1]}(x),\ \ldots,\ \tilde{u}^{[k-1]}(x),
\end{equation}
a $k$-step BDF applied to problem (\ref{eq:1}) can be written as
\begin{equation}\label{eq:bdf}
\Delta t b_k \mathcal{L}(\tilde{u}^{[n+k]}(x)) + a_k \tilde{u}^{[n+k]}(x)=
-\sum_{j=0}^{k-1}a_j \tilde{u}^{[n+j]}(x)
+ \Delta t b_k f(t_{n+k},x),
\end{equation}
\noindent
where $n=0,1,\ldots$ and the coefficient $a_j$, $j=0,\ldots,k$, and $b_k$ are listed in Table~\ref{tab:BDF} for $k \le 6$. It is worth noting that the 1-step BDF is the Backward Euler method.

BDF methods are $L$-stable for $k=1,2$ and $L(\alpha)$-stable for $k=3,\ldots,6$, with stability angle $\alpha$ reported in the last column of Table~\ref{tab:BDF}. It is also well known that they are not zero-stable for $k>6$. Definitions and further details can be found in~\cite{hairer1993solving}. The $L$-stability and $L(\alpha)$-stability make the methods suitable for ODEs whose solutions present high-frequency components. Furthermore, since BDF methods involve the right-hand side evaluation only at the right end of the current step, they allow the preservation of the elliptic structure in the time-discretized operator. In particular, the elliptic operator on the left-hand side is balanced by a source term involving known quantities that will be referred to as right-hand side and denoted by $\bm{RHS}^{[n]}$.

In order to give an estimate of the error in time, we follow \cite[Chapter~10]{thomeeLibro} and rewrite the BDF time discretization~\eqref{eq:bdf} of~\eqref{eq:1} as
\[\bar{\partial}_k\tilde u^{[n]}=\mathcal{L}(\tilde u^{[n]}) + f^{n},\]
with 
$\tilde u^{[n]} \approx u(n\Delta t,x)$, $f^{n}=f(n\Delta t,x)$. Here $\bar{\partial}_k$ is the backward difference operator given by $\displaystyle\bar{\partial}_k \tilde u^{[n]}=\frac{1}{\Delta t} \sum_{j=0}^k \frac{a_{k-j}}{b_k} \tilde u^{[n-j]}$. Letting $u^n=u(n\Delta t,x)$, we have
\[\bar{\partial}_k u^n=\mathcal L(u^n)+f^{n}+\tau^{n},\]
with $\tau^{n}=\bar{\partial}_k u^n-\partial_t u^n$ the consistency error. Then Lemma 10.1 and Theorem 10.1 in~\cite{thomeeLibro} provide the result reported next.
\begin{theorem}\label{th:error}
Let $k\leq 6$ and assume that the solution $u$ of~\eqref{eq:1} is sufficiently smooth. Then 
\[ \|u^n-\tilde u^{[n]}\|\leq c \sum_{j=0}^{k-1}\|u^j-\tilde u^{[j]}\|+ c \, \Delta t^k\int_{t_0}^{n\Delta t}\left\|\frac{\partial^{k+1}u}{\partial t^{k+1}}\right\|ds, \]
where $c$ is a positive constant.
\end{theorem}

\noindent
From Theorem \ref{th:error} it is clear that if the starting values~\eqref{eq:starting} are accurate enough, then the $k$-step BDF method~\eqref{eq:bdf} has order $k$.
Thus, a suitable strategy is needed to provide the first $k-1$ approximations. We refer to this strategy as \emph{starting procedure}. The starting values~\eqref{eq:starting} can be approximated by means of a $(k-1)$-step BDF method applied with a reduced stepsize $\Delta t / m$, with an integer, suitably chosen $m$.
As described in the previous sections, following~\cite{Cal_1} we choose the ELM collocation at given points for the discretization of the steady-state problems. The total error is then the time-marching error term plus the contribution of the ELM collocation error.

%
{\renewcommand{\arraystretch}{1.75}%
\begin{table}[t!]
\begin{center}
\caption{Coefficients and $L(\alpha)$-stability angles for the BDF methods.\label{tab:BDF}}
 \begin{tabular}{c|cccccccc||c}
 \hline
  $k$ & $a_6$ & $a_5$  & $a_4$ &  $a_3$ &  $a_2$ &  $a_1$ &  $a_0$ &  $b_k$ & $\alpha$ \\
 \hline
 1 &  &  &   &                  &                   &                 1 &              -1 & 1 & 90° \\
 2 &  &  &   &                  &                 1 &    $-\frac{4}{3}$ &   $\frac{1}{3}$ & $\frac{2}{3}$ & 90°\\
 3 &  &  &   &                1 &  $-\frac{18}{11}$ &    $\frac{9}{11}$ & $-\frac{2}{11}$ & $\frac{6}{11}$ & 86.03° \\
 4 &  &  & 1 & $-\frac{48}{25}$ &   $\frac{36}{25}$ &  $-\frac{16}{25}$ &  $\frac{3}{25}$ & $\frac{12}{25}$ & 73.35°\\
 5 &  & 1 & $-\frac{300}{137}$ & $\frac{300}{137}$ &   $-\frac{200}{137}$ &  $\frac{75}{137}$ &  $-\frac{12}{137}$ & $\frac{60}{137}$ & 51.84° \\
 6 & 1 & $-\frac{360}{147}$ & $\frac{450}{147}$ & $-\frac{400}{147}$ &   $\frac{225}{147}$ &  $-\frac{72}{147}$ &  $\frac{10}{147}$ & $\frac{60}{147}$ & 17.84°\\
 \hline
 \end{tabular}
\end{center}
\end{table}
}


%


In Algorithm~\ref{alg:elm-collocation} we describe the main steps required for the resolution. Since in our experiments we consider $d=1$, i.e.~a scalar variable $x$, we focus on this case. In principle the resolution at a fixed time step can be done with different methods, e.g.~one can change the number or the parameters of the involved AFs, or the number and the location of the collocation points. By the way, as can be seen from the pseudo-code (statements 2-4), if no changes are made, the linear problem that has to be solved at each time step involves the same matrix, so that this has to be constructed only once.

\begin{algorithm}[t]
\begin{algorithmic}
\caption{ELM collocation for the parabolic PDE with a {$k$-step BDF}\label{alg:elm-collocation}}
\State \textbf{Input:} the starting values \eqref{eq:starting}, the number of neurons $N$, the collocation points $x_j$, $j=1,\dots, M$. \;
\State \Comment{\textbf{Initialization}} \;
\State 1. Choose randomly $\alpha_i$ and $\beta_i$, $i=1,\ldots,N$, according to \eqref{eq:alpha}. \;
\State \Comment {\textbf{Discrete counterpart of the parabolic operator }} \;
\State 2. Compute the collocation matrix $C \in \R^{M\times N}$ by evaluating the left-hand side in~\eqref{eq:bdf} at the internal collocation points and appending the collocated boundary conditions. 
\State \Comment{\textbf{Time loop}}
\For{{$n = k:Nt$}}
\State 3. Compute $\bm{RHS}^{[n]} \in \R^{M} $ by evaluating the right-hand side in~\eqref{eq:bdf} at the collocation points. \;
\State 4. Find $\bm{w}^{[n]} \in \R^{N}$ that solves the linear problem $C \bm{w}^{[n]} = \bm{RHS}^{[n]}$. \;
\EndFor \;
\State
\State \textbf{Output:} the external ELM weights $w_i^{[n]}$, with $i = 1, \ldots, N$ and $n=k,\dots, Nt$, which provide the collocated solution $\tilde{u}^{[n]}(x)$ in~\eqref{eq:sol_meth_2}. \;
\end{algorithmic}
\end{algorithm}

\section{Numerical experiments\label{sec:experiments}}

Our test set consists of the following classes of problems:
\begin{enumerate}[label=(\alph*)]
\item stiff parabolic equations,
\item problems where the boundary conditions are discontinuous,
\item problems where the solution decays very rapidly.
\end{enumerate}
To obtain these problems, we used the one-dimensional heat equation with different diffusion coefficients and boundary conditions, which lead to challeging problems for the stability and accuracy properties of the methods~\cite{Cash}. The exact solution of these problems is known and we computed the approximate solution up to the final time $t_f$ using different values of the fixed time step length $\Delta t$. For the resolution at each time step, we used collocation with an ELM function consisting of $N$ neurons. Collocation was done on $M=N/2$ equispaced points by evaluating the exact derivatives of the activation functions, as in \cite{Cal_1}. For computing the final error, we used the $l_\infty$ norm of the difference between the approximate and the exact solution evaluated on $5000$ equispaced points in the spatial domain at the final time $t_f$.

All the numerical experiments were performed using MATLAB R2021b. The results were slightly affected by the choice of the parameters of the ELM functions, which were randomly generated as described in Section~\ref{sec:shallow-elms}. In our computations, these random parameters were obtained by using the MATLAB \texttt{randn} function. We found close behaviors with different initializations of \texttt{randn}, thus we decided to present the results obtained with \texttt{rng(1000)} for reproducibility issues, where \texttt{rgn} is the function specifying the seed for the random number generator.

In the plots we report the errors with respect to the number $Nt$ of stationary problems solved, which includes, for a fair comparison with one-step methods, the computational effort due to the starting procedures for the BDF methods. In the tests, we used the starting procedure described in Section~\ref{sec:time-marching}, with $m=8$. We chose to report the number of stationary problems solved because it represents the number of constructed underdetermined linear systems and hence it can be taken as representative of the overall computational cost. Our tests compare the behavior of Backward Euler (BE), Trapezoidal Rule (TR), and BDF of order 2-4 (referred to in the plots as BDF2-BDF4). We did not consider the higher-order methods BDF5 and BDF6 because if a modest number of collocation points is needed - as it is in our case - then the cost of the starting procedure dominates the overall performance. 

According with the discussion in Section~\ref{sec:shallow-elms}, the minimum-morm least squares solutions of the underdetermined linear systems were computed, by using the MATLAB \texttt{lsqminnorm} function with rank tolerance $10^{-15}$, which implements the COD algorithm. Of course, if the size of the collocation matrix $C$ is very large, this choice may be computationally very expensive, in terms of both time and memory. However, this is not the case of the experiments reported in this work, which are aimed at a methodological investigation.

\begin{figure}[h]%
\centering
\includegraphics[scale=0.62]{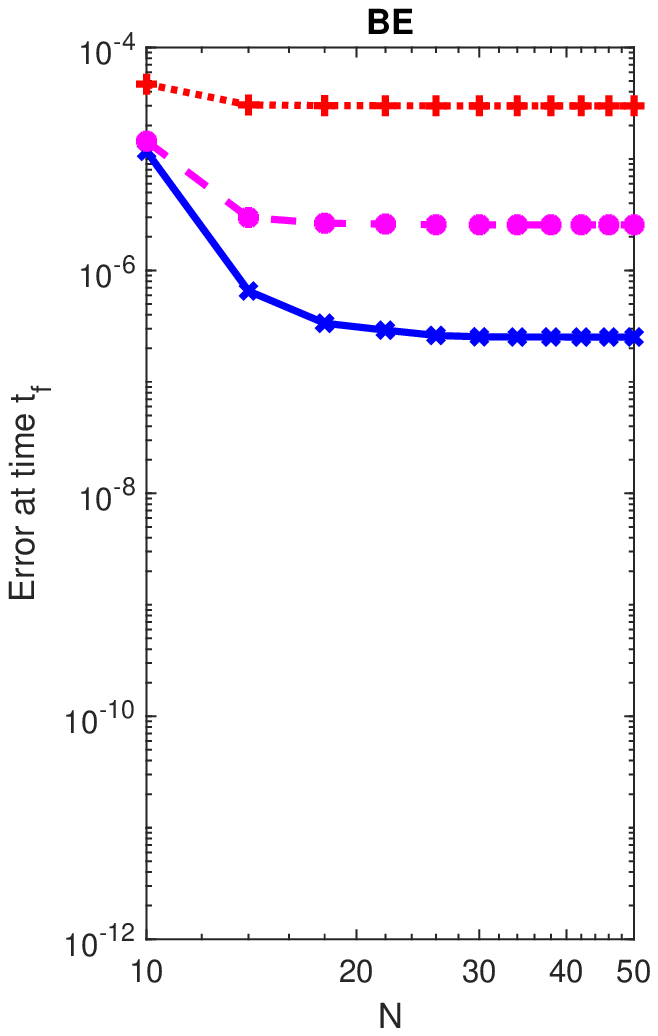}
\includegraphics[scale=0.62]{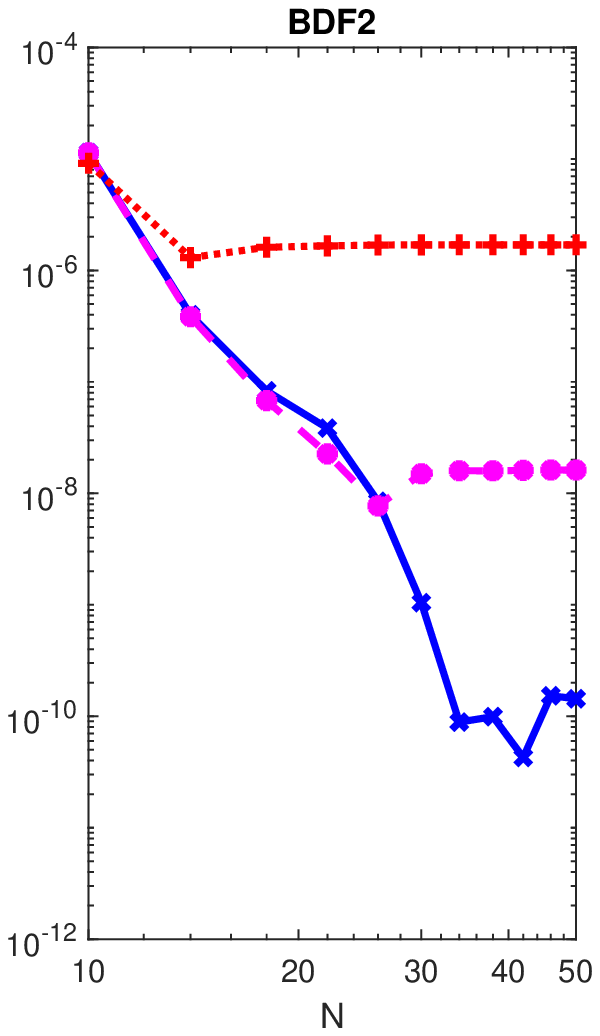}
\includegraphics[scale=0.62]{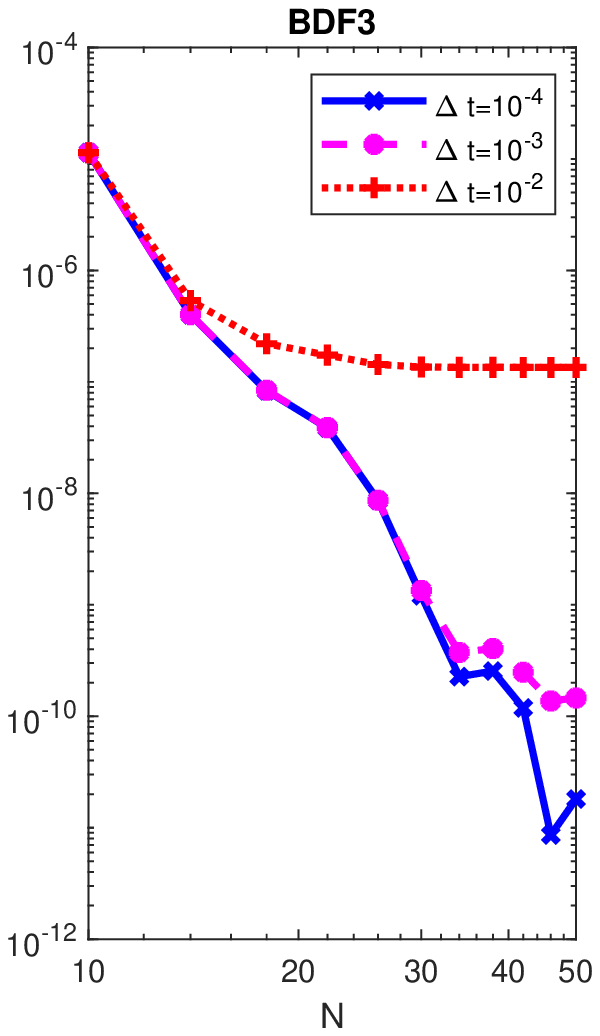}
\caption{Computed error for problem~\eqref{Pb-a} with $\gamma=3$ at the final time. We use the BE, BDF2 and BFD3 methods in time with tree different choices of $\Delta t$ and solve the problem by increasing the number of neurons $N$.\label{figNN}}
\end{figure}

\begin{figure}[ht]%
{\centering
\includegraphics[scale=0.62]{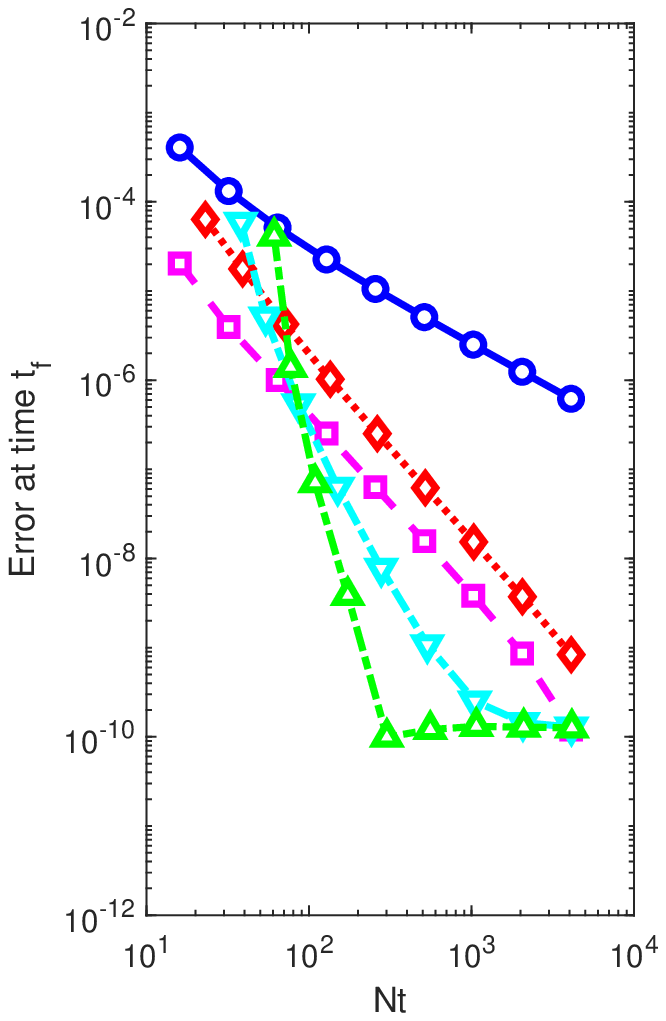}
\quad \quad
\includegraphics[scale=0.62]{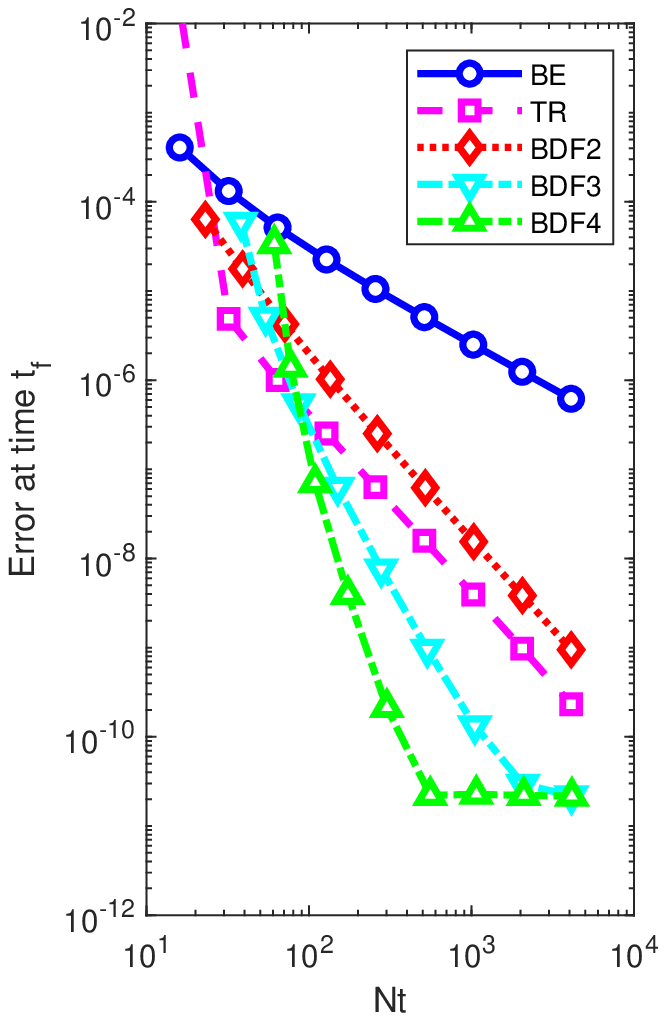}
\caption{Computed error for problem~\eqref{Pb-a} with $\gamma=3$ (left) and $\gamma=5$ (right) at the final time. The solutions are computed with $N=40$ neurons, decreasing $\Delta t$. In abscissae we report the number $Nt$ of linear systems solved, in ordinates the absolute value of exact error.\label{fig_reg}}}
\end{figure}



Test problem (a) is the following:
\begin{equation}\label{Pb-a}
\left\{
    \renewcommand{\arraystretch}{1.5}
    \begin{array}{ll}
       \dfrac{\partial u}{\partial t} = \dfrac{\partial^2 u}{\partial x^2}, & \ x\in[0,1], \ t \in[0,1], \\
        u(x,0) = \sin(\pi x) + \sin(\gamma \pi x), & \ u(0,t)=u(2,t)=0 ,
    \end{array} \right.
\end{equation}
whose exact solution, with $\gamma$ as a parameter, is:
\begin{equation*}
u(t,x)= e^{-\pi^2 t}\sin(\pi x) + e^{-\gamma^2 \pi^2 t}\sin(k \pi x) .
\end{equation*}
\noindent
When $\gamma$ increases, the second addend in the exact solution decays rapidly and oscillates. For this reason, the problem is referred to as stiff when $\gamma > 5$, while it is a ``standard'' test problem in the other cases.

The first numerical test that we present was carried out with $\gamma=3$ and regards the convergence of the method used to solve the stationary problems, i.e.~the ELM collocation when increasing the number $N$ of neurons. In Figure~\ref{figNN} we present the error computed at the final time $t_f$. In the different panels, we consider three different time-marching methods to solve the problem. Two facts can be noticed:
\begin{itemize}
    \item the convergence is very fast: we can conjecture a spectral convergence, compared also with the results obtained in~\cite{marcati2021exponential};
    \item the error stops decreasing when it reaches the maximum accuracy of the method in time: in that case, the latter begins to prevail.
\end{itemize} 
The three lines correspond to three different choices of $\Delta t$ and the fact that they coincide for small numbers of neurons confirms that at the beginning of the convergence history the error of the space discretizations prevails. The overall accuracy can be compared with that reported in the error plot in the left panel of Figure~\ref{fig_reg}, where the same problem is solved. {By looking also at the right panel of Figure~\ref{fig_reg}, corresponding to $\gamma=5$, we see that on these nonstiff problems our procedure has the expected behavior:} the order of convergence is achieved; the TR method slightly overperforms BDF2 because of its smaller error constant; the starting procedure for the BDF methods gives a shift of the initial points, which is negligible when $Nt$ grows; the error reaches its limit, given by the accuracy of the resolution by collocation.

\begin{figure}[ht]%
{\centering
\includegraphics[scale=0.62]{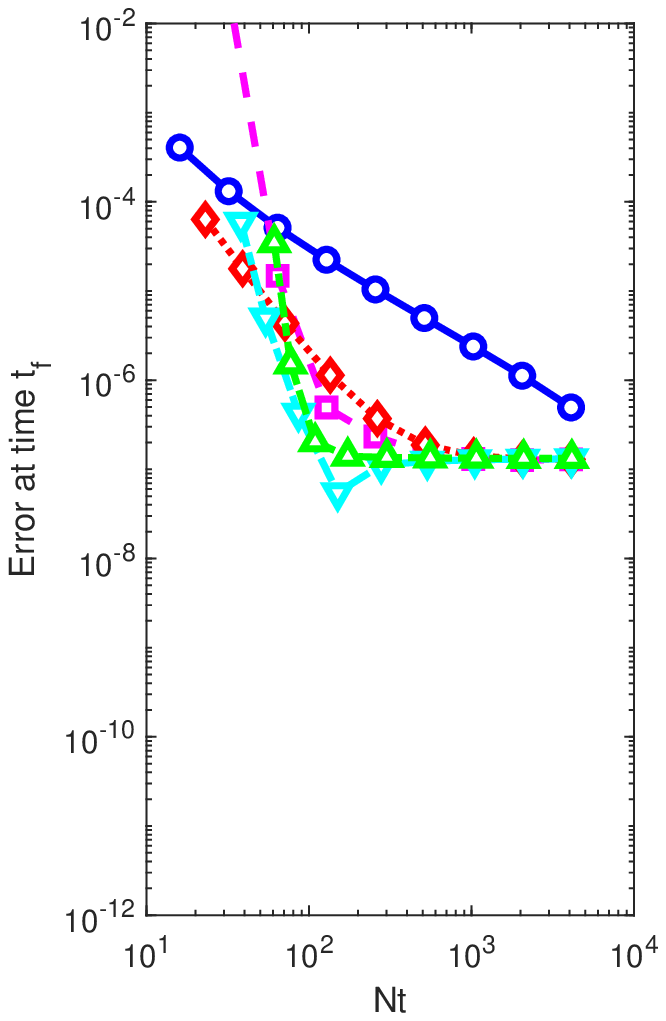}
\quad \quad
\includegraphics[scale=0.62]{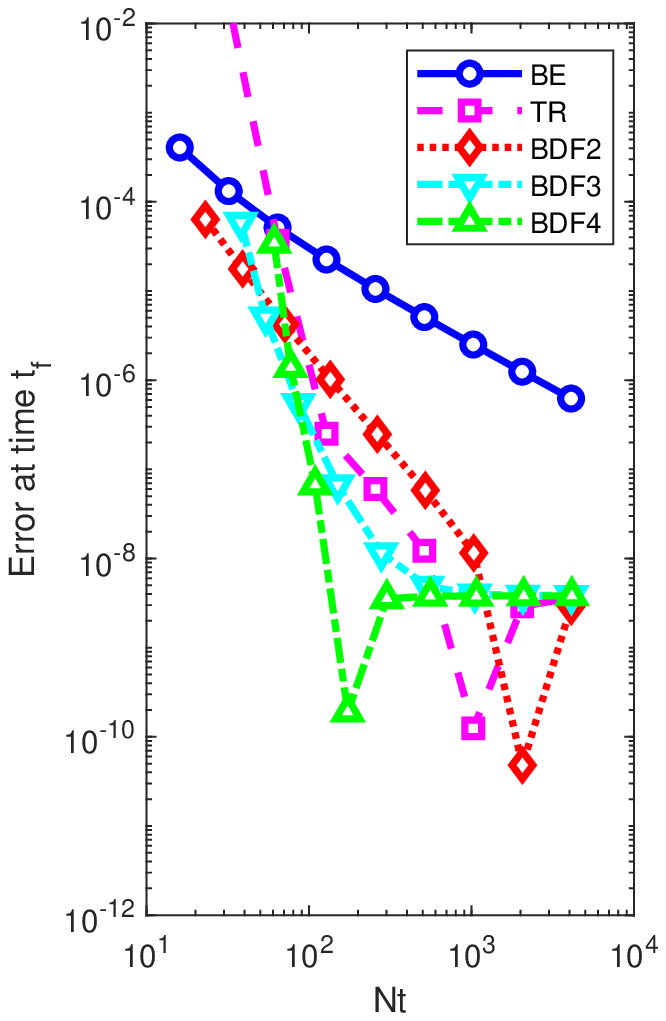}
\caption{Computed error for problem~\eqref{Pb-a} with $\gamma=10$ at the final time. The solutions are calculated with $N=40$ (left) and $N=50$ neurons (right), decreasing $\Delta t$. In abscissae we report the number $Nt$ of linear systems solved, in ordinates the absolute value of the error.\label{fig_stiff}}}
\end{figure}

In Figure~\ref{fig_stiff} we report the errors obtained while solving Problem~\eqref{Pb-a} with $\gamma=10$. The errors in the left panel correspond to $N=40$, while those in the right panel to $N=50$. This is a stiff problem, for which the trapezoidal method gives poor performances in the case of larger $\Delta t$ values, see also \cite{Cash}, while the BDF methods show a regular behavior because of their better stability properties. This is also a difficult problem for the collocation method, being the exact solution highly oscillating. For this reason the final accuracy is poor for $N=40$ and improves for $N=50$. 

Problem (b) (with discontinuous boundary conditions) is:
\begin{equation}\label{Pb-b}
\left\{ \renewcommand{\arraystretch}{1.5}
    \begin{array}{ll}
       \dfrac{\partial u}{\partial t} = \dfrac{\partial^2 u}{\partial x^2},  & \ x\in[0,2], \ t \in[0,1.2], \\
        u(x,0)=1, & \ u(0,t)=u(2,t)=0.
    \end{array}\right.
\end{equation}
\noindent
The exact solution\footnote{The solution is given in the form of a series, but its terms decay very rapidly, so that for our numerical tests the approximate solution obtained with 20 terms is exact up to machine precision.} to problem~\eqref{Pb-b} is:
\begin{equation*}
u(t,x)=\sum_{n=1}^{\infty} \left[1-(-1)^n \right] \dfrac{2}{n\pi} \sin\left( \dfrac{n\pi x}{2} \right) \exp\left( \dfrac{-n^2\pi^2 t}{4} \right) .
\end{equation*}

\begin{figure}[th]%
{\centering
\includegraphics[scale=0.62]{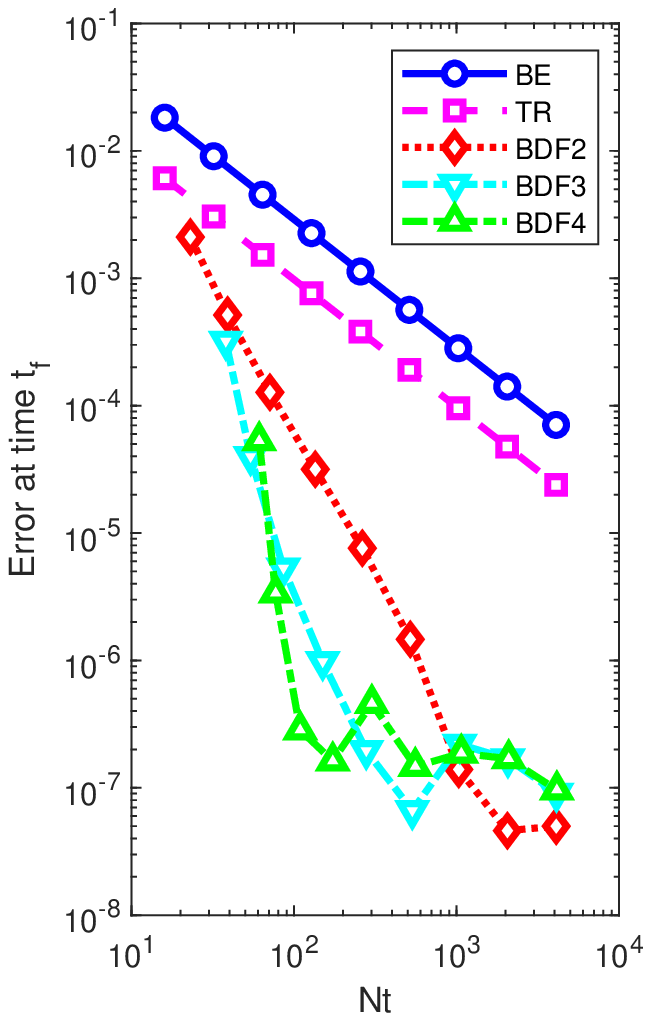}
\quad \quad
\includegraphics[scale=0.62]{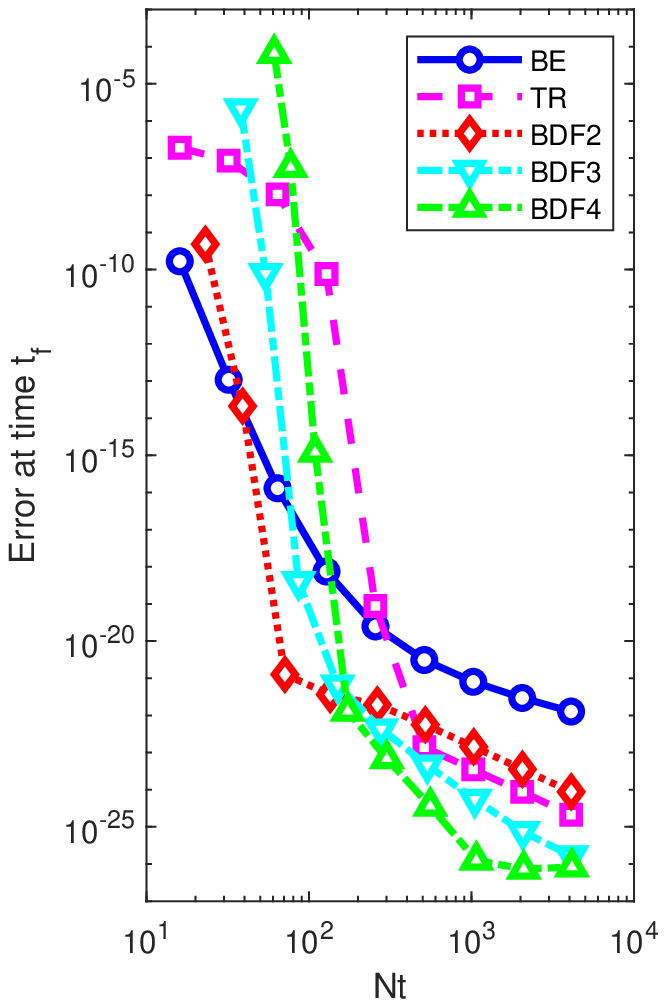}
\caption{Computed error for problem~\eqref{Pb-b} (left panel) and for problem~\eqref{Pb-c} (right panel) at the final time. The solutions are calculated with $N=40$ neurons, decreasing $\Delta t$. In abscissae we report the number $Nt$ of linear systems solved, in ordinates the absolute value of the error.\label{fig_2test_finali}}}
\end{figure}

Numerical results for this test case are reported in the left panel of Figure~\ref{fig_2test_finali}. As already noticed in~\cite{Cash}, the trapezoidal rule suffers from order reduction, behaving as an order-1 method, while the BDF solvers behave as expected. Moreover, the BDF accuracy in this case is much higher with the same computational cost.

Problem (c) (with a solution that decays very rapidly) is:
\begin{equation}\label{Pb-c}
    \left\{ \renewcommand{\arraystretch}{1.5}
    \begin{array}{ll}
       \dfrac{\partial u}{\partial t} = 5 \dfrac{\partial^2 u}{\partial x^2}, & \ x\in[0,1], \ t \in[0,1], \\
        u(x,0)=\sin(\pi x), & \ u(0,t)=u(1,t)=0,
    \end{array}\right.
\end{equation}
and its solution is:
\begin{equation*}
u(t,x)= e^{-5\pi^2 t}\sin(\pi x).
\end{equation*}

Numerical results for this test case are reported in Figure~\ref{fig_2test_finali}, in the right panel. One can observe once again that the methods can achieve very high accuracy and maintain the properties seen in the standard test cases.

\section{Conclusions and future work\label{sec:conclusions}}

Scientific Machine Learning (SML) is a research field in which Artificial Intelligence methodologies have been employed to solve in innovative manners problems modeled by PDEs. A large number of easy-to-use methods, based on ANNs, allow researchers to deal with complex PDEs efficiently. In this paper, we designed a numerical scheme belonging to the class of ELM methods for solving time-dependent parabolic PDE problems. This can be considered as a step towards rethinking meshless methods by using ANNs. The proposed approach addresses the time evolution by applying time-marching techniques and adopts the collocation for solving the resulting stationary elliptic problems. Regarding the elliptic component of the PDE, we observed a good accuracy with a limited computational costs in situations where classical methods fail to give good results.

The main goal of this work was to explore how the time approximation affects the resolution of time-dependent PDEs by ELM-based collocation, using classical methods to discretize the problem in time. To this aim, some linear parabolic PDE that are challeging problems for the stability and accuracy properties of the methods 
were considered. The $\theta$-method and BDF techniques were investigated. We observed that BDF methods have good accuracy and convergence properties while keeping the same computational cost as the implicit $\theta$-method. Moreover, the time-discretization schemes used in the space-collocation method have promising properties in terms of stability and order of convergence, which are confirmed by numerical tests. We also concluded that, in our numerical framework, BDF methods of order 2-4 were to be preferred to Backward Euler or Trapezoidal rules, because the latter are unable to deal with high accuracy requirements, as they are slowly convergent and can suffer from order reduction or numerical instability.

SML methodologies to design novel numerical methods for solving PDEs represent a fascinating research field for which contributions grow exponentially. Among future improvements of the proposed approach we identified two main directions: i) network architectural studies and ii) theoretical results. Concerning direction i), how integrating ELMs with physics-informed approaches like PINNs is a challenging task; concerning ii), the adoption of ANNs as universal approximators of PDE solutions has to be supported by theoretical results about errors in the learning process, such as a-priori bounds related to stability and convergence rate.



\section*{Declarations}

\bmhead{Funding}
This work was partially supported by the Istituto Nazionale di Alta Matematica - Gruppo Nazionale per il Calcolo Scientifico (INdAM-GNCS), Italy.

\bmhead{Data availability}
Data sharing is not applicable to this article as no datasets were generated during the current study.

\bmhead{Conflict of interest}
The authors declare that they have no conflict of interest.



\end{document}